\documentclass{article}
\usepackage{amsmath,amsthm}
\usepackage{graphicx,epstopdf,epsfig,multirow,epic,bm}
\usepackage{lineno,hyperref}
\usepackage{amssymb}

\normalsize \rm
\begin{document}

\begin{center}
{\Large  \textbf {On the existence of neutral graph}}\\[12pt]
{\large Fei Ma$^{a,}$\footnote{~The author's E-mail: feima@nwpu.edu.cn. } }\\[6pt]
{\footnotesize $^{a}$  School of Computer Science, Northwestern Polytechnical University, Xi'an 710072, China}\\[12pt]
\end{center}

\begin{quote}
\textbf{Abstract:} Graph is considered neutral if its assortativity coefficient $r$ is equal to zero. In this paper, we address an outstanding conjecture, i.e., whether is there a neutral graph on $n$ vertices? First, we show that for $n\geq7$, there is at least one neutral tree, which suggests that we find a representative of any order neutral graph. Additionally, we obtain that given $n\geq13$, there exists at least one neutral non-tree graph.\\

\textbf{Keywords:} Neutral graph, Assortativity coefficient. \\

\end{quote}

\vskip 1cm

\section{Introduction}

Let $G=(V,E)$ represent a simple connected graph, then, notations $V$ and $E$, respectively, indicates vertex set and edge set. We use $G$ as shorthand of $G=(V,E)$ when it is clear from the context. The number of vertices (a.k.a., order) of graph $G$ is $n(:=|V|)$, and the total number of edges in graph $G$ is $m(:=|E|)$. Assorataivity coefficient $r$ of graph $G$ is defined as follows \cite{Newman-2002}

\begin{equation}\label{eqa:MF-JOCT-0}
r=\frac{m^{-1}\sum\limits_{e_{uv}\in E} d_{u}d_{v}-\left[m^{-1}\sum\limits_{e_{uv}\in E} \frac{1}{2}(d_{u}+d_{v})\right]^{2}}{m^{-1}\sum\limits_{e_{uv}\in E} \frac{1}{2}(d^{2}_{u}+d^{2}_{v})-\left[m^{-1}\sum\limits_{e_{uv}\in E} \frac{1}{2}(d_{u}+d_{v})\right]^{2}},
\end{equation}
in which $d_{u}$ denotes degree of vertex $u$ and $e_{uv}$ represents an edge whose endpoints are vertices $u$ and $v$, separately. Graph $G$ is considered assortative if $r>0$, disassortative if $r<0$, and neutral otherwise. Obviously, it is easy to determine whether a given graph is neutral or not. However, to the best of our knowledge, there have not been published work answering the following conjecture \cite{Noldus-2015}-\cite{Behrens-2022}
$$\text{Is there neutral graph of arbitrary order?}$$

In this work, we provide a positive answer to the aforementioned issue.

\section{Definition and main results}

This section comprises two subsections as follows.

\subsection{Definition} First of all, we introduce some notations and terminologies used throughout in this work. Note that graph under consideration is simple and connected.

\textbf{Definition 1}\cite{Bondy-2008} A graph $G$ is considered $k$-regular if and only if each vertex $u$ has degree $d_{u}=k$. For simplicity, we denote by $k$-$G$ $k$-regular graph. $k$-$G$ is a cycle if $k$ is equal to $2$.

\textbf{Definition 2}\cite{Bondy-2008} The \emph{$s$-th subdivision graph} $G_{s}=(V_{s}, E_{s})$ is obtained from graph $G$ by inserting $s$ vertices into each edge.

\textbf{Definition 3} The \emph{single-edge division graph} $G'=(V', E')$ is obtained from graph $G$ by only inserting one vertex into an arbitrarily chose edge.

\textbf{Definition 4} The \emph{triangle-operation} is to (i) first add an isolated vertex for each edge of graph $G$, and (2) connect each isolated vertex to two endvertices of the associated edge using two edges. The resulting graph is referred to as $G^{\bigtriangleup}=(V^{\bigtriangleup}, E^{\bigtriangleup})$ for brevity.

\textbf{Definition 5} The \emph{leaf-connecting-operation} is to (1) connect $d_{v}$ leaf vertices to each vertex $v$ in graph $G$, and (2) connect a pair of leaves using an edge, then perform the same connection between another two leaves until no leaf is left. The resulting graph is referred to as $G^{\dagger}=(V^{\dagger}, E^{\dagger})$ for presentation. Degree one vertex is defined as a leaf.

\subsection{Main results} In the sequel, let us show main results composed of two subsections. The former shows the existence of neutral tree, and the existence of neutral non-tree graph is displayed in the later.

\subsubsection{Neutral tree} Here we prove the existence of neutral tree.

\textbf{Lemma 1} Graph $G$ is either a cycle of any order or a tree with a unique largest degree $\Delta=3$ vertex, if the following expression holds

\begin{equation}\label{eqa:MF-JOCT-1-1}
\sum_{e_{uv}\in E}(d_{u}+ d_{v})=4m.
\end{equation}
From which it follows that the desired tree has order $n\geq4$.

\emph{Proof} It is easy to prove Lemma 1 using some elementary computations. \qed

\textbf{Lemma 2} Given an arbitrary graph $G$, the following expression holds 

\begin{equation}\label{eqa:MF-JOCT-1-2}
2m\sum\limits_{e_{uv}\in E}(d^{2}_{u}+d^{2}_{v})-\left[\sum\limits_{e_{uv}\in E}(d_{u}+d_{v})\right]^{2}\geq0.
\end{equation}
Equality holds iff $G$ is regular.

\emph{Proof} By means of Cauchy-Schwartz inequality, we can easily prove Lemma 2. \qed

Based on Lemma 2, it is clear to see that an irregular graph $G$ is neutral iff the numerator of Eq. (\ref{eqa:MF-JOCT-0}) is equal to zero.

\textbf{Lemma 3} In $s$-th subdivision graph $G_{s}=(V_{s}, E_{s})$,

\begin{equation}\label{eqa:MF-JOCT-1-3}
\begin{aligned}4|E_{s}|\sum\limits_{e_{ij}\in E_{s}} d_{i}d_{j}-\left[\sum\limits_{e_{ij}\in E_{s}}(d_{i}+d_{j})\right]^{2}=-\left[\sum_{e_{uv}\in E}(d_{u}+ d_{v}-4)\right]^{2}.
\end{aligned}
\end{equation}

\emph{Proof} Using some simple algebra, we complete Lemma 3. \qed

\textbf{Theorem 1} Given any tree $T$ subjected to Eq.(\ref{eqa:MF-JOCT-1-1}), the corresponding $s$-th subdivision tree $T_{s}$ is neutral.

\emph{Proof} Using Lemmas 1 and 3, we verify Theorem 1. \qed

According to Theorem 1, together with Lemma 1, any tree of order $2n-1$ turns out to be neutral by setting $s=1$. In addition, using simple calculations, we reach the following corollary.

\textbf{Corollary 1} Any tree $T'_{s}$ is neutral, and its single-edge division tree $(T'_{s})'$ is neutral as well.

Now, based on Theorem 1 and Corollary 1, it is easy to prove that there is at least one neutral tree of order $n+3$. Additionally, we easily verify, by enumeration,  that any graph is not neutral if its order is no more than $6$. Armed with the results above, we firmly state that the next theorem holds true.

\textbf{Theorem 2} There exists at least one tree of order $n+3$ such that the associated assortativity coefficient $r$ is equivalent to zero.

By far, we succeed in finding a representative of any order neutral graph. It should be mentioned that the representative has tree structure. Next, we are concerned with other representatives having non-tree structure.

\subsubsection{Neutral non-tree graph} Below elaborates on the existence of neutral non-tree graph.

\textbf{Lemma 4} Given an arbitrary neutral tree $T$ on $n_{1}(\geq7)$ vertices, then the graph $T^{\dagger}$ is neutral.

\emph{Proof} After some simple algebra, we complete Lemma 4. \qed

Clearly, neutral graph $T^{\dagger}$ is not a tree. It should be mentioned that graph $T^{\dagger}$ has $3n_{1}-2$ vertices.

\textbf{Lemma 5} Given an arbitrary cycle $G$ on $n_{2}(\geq3)$ vertices, then the graph $G^{\dagger}$ is neutral.

\emph{Proof} As before, with some simple algebra, we complete Lemma 5. \qed

Similarly, neutral graph $G^{\dagger}$ is not a tree either. At the same time, graph $G^{\dagger}$ has $3n_{2}$ vertices.

\textbf{Lemma 6} Given an arbitrary neutral tree $T$ on $n_{1}(\geq7)$ vertices, then the graph $T^{\bigtriangleup}$ is neutral.

\emph{Proof} Analogously, with some simple algebra, we complete Lemma 6. \qed

It is clear to see that non-tree graph $T^{\bigtriangleup}$ has $2n_{1}-1$ vertices.

Until now, we obtain that there is a neutral non-tree graph where the total number of vertices $n$ obeys

\begin{equation}\label{eqa:MF-JOCT-2-1}
n\mod 6=0,1,3,4,5, \qquad n\geq\max\left\{\min\{3n_{1}-2\},\min\{3n_{2}\},\min\{2n_{1}-1\}\right\}=19
\end{equation}

Next, we will prove case $n\mod6=2$. To this end, we need to introduce a class of graphs. Without loss of generality, we use $G^{\oplus}$ to represent these graphs. Graph $G^{\oplus}$ is obtained based on two small graphs. Concretely speaking, given a neutral graph $G(1)=(V(1),E(1))$ and a $k$-regular graph $G(2)=(V(2),E(2))$, then we require that (i) $\beta$ incomplete edges\footnote{An edge is considered incomplete if it just connects a vertex.} be stemmed from each vertex of graph $G(2)$ and (ii) for each edge of graph $G(1)$, $(\alpha-1)$ incomplete edges be stemmed from every endvertex. In the meantime, we assume that $2(\alpha-1)|E(1)|=\beta|V(2)|$ holds true. Next, one merges an incomplete edge from graph $G(1)$ with another incomplete edge in graph $G(2)$ into an edge. Performing such manipulation until no incomplete edge is left. The resultant graph is denoted by $G^{\oplus}$.

\textbf{Lemma 7} Graph $G^{\oplus}$ is neutral if the following expression holds

\begin{equation}\label{eqa:MF-JOCT-2-2}
\frac{|E(2)|}{|E(1)|}=(\alpha-1)^{2}=\left(\frac{k}{\beta}\right)^{2}=\left(\frac{|V(2)|}{|V(1)|}\right)^{2}.
\end{equation}

\emph{Proof} By definition in (\ref{eqa:MF-JOCT-0}), we have

\begin{equation}\label{eqa:MF-JOCT-2-2-1}
\begin{aligned}
\sum\limits_{e_{ij}\in E^{\oplus}} d_{i}d_{j}&=\sum\limits_{e_{uv}\in E(1)} \alpha^{2}d_{u}d_{v}+\sum\limits_{u\in V(1)} \alpha d_{u}(k+\beta)(\alpha-1)d_{u}+(k+\beta)^{2}|E(2)|\\
&=\sum\limits_{e_{uv}\in E(1)} \alpha^{2}d_{u}d_{v}+\alpha(\alpha-1)(k+\beta)\sum\limits_{u\in V(1)}d_{u}^{2}+(k+\beta)^{2}|E(2)|.
\end{aligned}
\end{equation}

\begin{equation}\label{eqa:MF-JOCT-2-2-2}
\begin{aligned}
\sum\limits_{e_{ij}\in E^{\oplus}} d_{i}+d_{j}&=\sum\limits_{e_{uv}\in E(1)} (\alpha d_{u}+\alpha d_{v})+\sum\limits_{u\in V(1)} [\alpha d_{u}+(k+\beta)](\alpha-1)d_{u}\\
&\quad+2(k+\beta)|E(2)|\\
&=\alpha^{2}\sum\limits_{e_{uv}\in E(1)}(d_{u}+d_{v})+2(k+\beta)[(\alpha-1)|E(1)|+|E(2)|].
\end{aligned}
\end{equation}

Next, it is easy to check

\begin{equation}\label{eqa:MF-JOCT-2-2-3}
\begin{aligned}
&4|E^{\oplus}|\left(\sum\limits_{e_{ij}\in E^{\oplus}} d_{i}d_{j}\right)-\left(\sum\limits_{e_{ij}\in E^{\oplus}} d_{i}+d_{j}\right)^{2}\\
&=\alpha^{4}\left[4|E(1)|\left(\sum\limits_{e_{uv}\in E(1)} d_{u}d_{v}\right)-\left(\sum\limits_{e_{uv}\in E(1)} d_{u}+d_{v}\right)^{2}\right]\\
&\equiv 0,
\end{aligned}
\end{equation}
in which $|E^{\oplus}|=(2\alpha-1)|E(1)|+|E(2)|$. This is complete. \qed

From the results above, we obtain

$$|V^{\oplus}|=|V(1)|+|V(2)|=\alpha|V(1)|.$$

After that, if we suppose that parameter $\alpha$ is equal to $2$, and neutral graph $G(1)$ is graph $T^{\dagger}$, the resulting graph $G^{\oplus}$ satisfies

$$|V^{\oplus}|=2(3n_{1}-2)\mod 6=2.$$

To sum up, we come to the following theorem.

\textbf{Theorem 3} Given $n\geq38$, there exists at least one neutral non-tree graph of order $n$.

As a byproduct, we have the following assertion about graph $G^{\ominus}$. The $G^{\ominus}$ is obtained from an arbitrary neutral graph $G$ and its copy $G^{\circ}$ in a similar manner as used in the construction of graph $G^{\oplus}$. Concretely speaking, one connects each vertex $u$ in graph $G$ to neighboring vertices of its peer $u^{\circ}$ in copy $G^{\circ}$.

\textbf{Assertion 1} Graph $G^{\ominus}$ shares the same assortativity coefficient with graph $G$.

\emph{Proof} By definition in (\ref{eqa:MF-JOCT-0}), we have

\begin{equation}\label{eqa:MF-JOCT-2-3-1}
\begin{aligned}
\sum\limits_{e_{ij}\in E^{\ominus}} d_{i}d_{j}&=4\sum\limits_{e_{uv}\in E} 2^{2}d_{u}d_{v}=16\sum\limits_{e_{uv}\in E} d_{u}d_{v},
\end{aligned}
\end{equation}

\begin{equation}\label{eqa:MF-JOCT-2-3-2}
\begin{aligned}
\sum\limits_{e_{ij}\in E^{\ominus}} d_{i}+d_{j}&=4\sum\limits_{e_{uv}\in E} (2 d_{u}+2d_{v})=8\sum\limits_{e_{uv}\in E} (d_{u}+d_{v}),
\end{aligned}
\end{equation}

\begin{equation}\label{eqa:MF-JOCT-2-3-3}
\begin{aligned}
\sum\limits_{e_{ij}\in E^{\ominus}} d^{2}_{i}+d^{2}_{j}&=4\sum\limits_{e_{uv}\in E} [(2d_{u})^{2}+(2d_{v})^{2}]=16\sum\limits_{e_{uv}\in E} (d^{2}_{u}+d^{2}_{v}),
\end{aligned}
\end{equation}

Next, it is easy to check

\begin{equation}\label{eqa:MF-JOCT-2-3-4}
\begin{aligned}
4|E^{\ominus}|&\left(\sum\limits_{e_{ij}\in E^{\ominus}} d_{i}d_{j}\right)-\left(\sum\limits_{e_{ij}\in E^{\ominus}} d_{i}+d_{j}\right)^{2}\\
&=64\left[4|E|\left(\sum\limits_{e_{uv}\in E} d_{u}d_{v}\right)-\left(\sum\limits_{e_{uv}\in E} d_{u}+d_{v}\right)^{2}\right],
\end{aligned}
\end{equation}
and

\begin{equation}\label{eqa:MF-JOCT-2-3-5}
\begin{aligned}
2|E^{\ominus}|&\left[\sum\limits_{e_{ij}\in E^{\ominus}} (d^{2}_{i}+d^{2}_{j})\right]-\left(\sum\limits_{e_{ij}\in E^{\ominus}} d_{i}+d_{j}\right)^{2}\\
&=64\left[2|E|\sum\limits_{e_{uv}\in E} (d^{2}_{u}+d^{2}_{v})-\left(\sum\limits_{e_{uv}\in E} d_{u}+d_{v}\right)^{2}\right].
\end{aligned}
\end{equation}
in which we have used quality $|E^{\ominus}|=4|E|$. By far, it is easy to prove $r^{\ominus}\equiv r$. This is complete. \qed

Obviously, if we suppose that graph $G$ is graph $T^{\dagger}$, the resulting graph $G^{\ominus}$ also satisfies

$$|V^{\ominus}|=2(3n_{1}-2)\mod 6=2.$$

In addition, if an arbitrary tree $T$ on $n_{1}(\geq7)$ vertices stated in Theorem 2 is selected as graph to create graph $G^{\ominus}$, then the total number of vertices of the end graph $G^{\ominus}$ is equal $2n_{1}$. Together with Lemma 6, we reach the last theorem in this work as follows.

\textbf{Theorem 4} Given $n\geq13$, there exists at least one neutral non-tree graph of order $n$.

\section*{Acknowledgments}
The research was supported in part by the National Natural Science Foundation of China No. 62403381, the Key Research and Development Plan of Shaanxi Province No. 2024GX-YBXM-021, and the Fundamental Research Funds for the Central Universities No. G2025KY06223

\end{document}